\documentclass{article}
\usepackage{amsmath}
\usepackage{amsfonts}
\usepackage{amssymb}
\usepackage{pstricks,pst-node,amssymb,amsmath,graphics,latexsym,tabularx,shapepar}
\usepackage{graphicx}
\usepackage{float}
\newtheorem{theorem}{Theorem}

\newcommand{\bpartial}{\mathop{\partial\kern -4pt\raisebox{.8pt}{$|$}}}
\newcommand{\bra}{\mathopen{[\kern-1.6pt[}}
\newcommand{\ket}{\mathclose{]\kern-1.5pt]}}
\newcommand{\bbra}{\mathopen{[\kern-2.2pt[\kern-2.3pt[}}
\newcommand{\bket}{\mathclose{]\kern-2.1pt]\kern-2.3pt]}}

\makeindex

\begin{document}

\title{Lie symmetry analysis of a (2+1)-dimensional flux-limited Keller--Segel system}
\author {\small{ \bf  Ahmed Abbas Jaber Al-Furaiji, \hspace{1mm} Ghorbanali Haghighatdoost\thanks{Corresponding author, \em gorbanali@azaruniv.ac.ir}\; and \hspace{1mm} Mustafa Bazghandi}\\
{\small{\em Department of Mathematics,Azarbaijan Shahid Madani University, Tabriz, Iran}}\\
{\thanks{ E-mails: \em ahmedgon1234567@gmail.com, \em mostafabazghandi2001@gmail.com}}
}
\maketitle

\begin{abstract}
We investigate a two-dimensional flux-limited Keller--Segel (FLKS) system using classical Lie symmetry analysis. The complete Lie point symmetry algebra is determined and shown to consist of temporal and spatial translations together with planar rotations. An optimal system of one-dimensional subalgebras is constructed, yielding stationary, translationally invariant, radially symmetric, travelling-wave, and rotating-wave reductions. The corresponding reduced equations are derived systematically. The results demonstrate that the flux-limiting mechanism eliminates the scaling symmetries of the classical Keller--Segel model and significantly restricts the class of admissible invariant solutions. Exact spatially homogeneous equilibrium solutions are determined, and integral representations for stationary radial and travelling-wave profiles are derived.

\end{abstract}

\noindent{\bf AMS}:35B06, 76M60, 35Q92.

\section{Introduction}

Chemotaxis, the directed movement of biological organisms in response to chemical stimuli, plays a fundamental role in numerous biological processes, including bacterial aggregation, embryonic development, tumor invasion, and immune responses. Since its introduction by Keller and Segel, the Keller--Segel (KS) model has become one of the most influential mathematical frameworks for describing chemotactic behavior through coupled reaction--diffusion equations governing the evolution of cell density and chemoattractant concentration \cite{KellerSegel1970}. Despite its success in capturing aggregation phenomena, the classical KS model is known to exhibit finite-time blow-up and unbounded cell accumulation in certain parameter regimes, particularly in higher spatial dimensions \cite{Buseghin2026}.

To overcome these limitations and provide a more realistic description of biological transport, several modified chemotaxis models have been proposed. Among them, flux-limited Keller--Segel (FLKS) systems have attracted considerable attention because they incorporate bounded chemotactic velocities through nonlinear flux-limiting mechanisms \cite{Furaiji2026,Marras2025,Mao2026}. Such modifications prevent unrealistic transport speeds and lead to richer dynamical behavior while preserving essential biological features. The resulting equations, however, possess a significantly more complicated nonlinear structure, making the construction of analytical solutions and the investigation of qualitative properties considerably more challenging.

Symmetry methods provide powerful tools for studying nonlinear differential equations. In particular, Lie group analysis enables the systematic determination of continuous transformation groups admitted by a given model and facilitates the construction of invariant solutions, similarity reductions, conservation laws, and exact analytical structures (e.g. \cite{Navratil2026,Samreen2026,Solberg2026,Yadav2026}). The fundamental theory of Lie symmetry methods and its application to PDEs are comprehensively presented in the monographs \cite{bluman2010,olver1993}. Recent developments have further demonstrated that symmetry-based approaches can reveal nonlocal structures that are not immediately apparent from the classical Lie symmetries (e.g. \cite{Bazghandi2026,BlumanYang2013}). 

In this paper, we consider the two-dimensional flux-limited Keller--Segel system

\begin{equation}
\begin{aligned}
u_t &= \nabla\cdot\left( \nabla u - \frac{\chi\,u\,\nabla v}{\sqrt{1+|\nabla v|^2}} \right), \\
v_t &= \Delta v + u - v,
\end{aligned}
\label{eq:FLKS}
\end{equation}
where $\chi>0$ denotes the chemotactic sensitivity, $(x,y,t)\in \mathbb{R}^2\times\mathbb{R}^+$, and $u(x,y,t)$ and $v(x,y,t)$ denote cell and chemoattractant concentrations, respectively.
 The flux limiter
$
\Phi(\nabla v)=\nabla v/\sqrt{1+|\nabla v|^2}
$
ensures bounded transport velocities and preserves rotational invariance, while simultaneously breaking the classical scaling symmetries of the standard KS model.


The main objective of this work is to perform a complete Lie point symmetry analysis of the FLKS system \eqref{eq:FLKS}. We determine the admitted Lie symmetry algebra, construct an optimal system of one-dimensional subalgebras, and derive the corresponding similarity reductions. These reductions yield stationary, translationally invariant, radially symmetric, travelling-wave, and rotating-wave formulations of the original model. In addition, we identify exact invariant solutions and discuss the influence of the flux-limiting mechanism on the symmetry properties of the system.

The paper is organized as follows. In Section~\ref{sec:symmetries}, we obtain the Lie point symmetries of the FLKS system \eqref{eq:FLKS}. Section~\ref{sec:optimal} is devoted to the construction of the optimal system of one-dimensional subalgebras. The resulting symmetry reductions and reduced equations are presented in Section~\ref{sec:reductions}. Section~\ref{sec:solutions} discusses invariant solutions of the reduced equations. Finally, Section~\ref{sec:conclusion} summarizes the main findings and outlines directions for future work.
\section{Lie Point Symmetries}\label{sec:symmetries}

To determine the Lie point symmetries of system \eqref{eq:FLKS}, we introduce the infinitesimal generator

\[
\mathbf{X} = \xi^1(x,y,t,u,v)\partial_x + \xi^2(x,y,t,u,v)\partial_y + \tau(x,y,t,u,v)\partial_t + \eta^1(x,y,t,u,v)\partial_u + \eta^2(x,y,t,u,v)\partial_v.
\tag{2}
\]

The second prolongation $\mathrm{pr}^{(2)}\mathbf{X}$ is given by

\[
\mathrm{pr}^{(2)}\mathbf{X} = \mathbf{X} + \eta^{1,x}\partial_{u_x} + \eta^{1,y}\partial_{u_y} + \eta^{1,t}\partial_{u_t} + \eta^{2,x}\partial_{v_x} + \eta^{2,y}\partial_{v_y} + \eta^{2,t}\partial_{v_t} + \cdots,
\]

where the prolonged coefficients follow the standard formulas

\[
\begin{aligned}
\eta^{1,x} &= D_x\eta^1 - (D_x\xi^1)u_x - (D_x\xi^2)u_y - (D_x\tau)u_t,\\
\eta^{1,y} &= D_y\eta^1 - (D_y\xi^1)u_x - (D_y\xi^2)u_y - (D_y\tau)u_t,\\
\eta^{1,t} &= D_t\eta^1 - (D_t\xi^1)u_x - (D_t\xi^2)u_y - (D_t\tau)u_t,
\end{aligned}
\]

and analogously for $\eta^{2,x},\eta^{2,y},\eta^{2,t}$. Here $D_x,D_y,D_t$ denote total derivatives.

The invariance conditions are

\[
\left.\mathrm{pr}^{(2)}\mathbf{X}(\Delta_i)\right|_{\Delta_1=0,\;\Delta_2=0} = 0,\qquad i=1,2.
\tag{5}
\]

\subsubsection*{Step 1: Determining equations from $\Delta_2$}

The second equation $\Delta_2 = v_t - v_{xx} - v_{yy} - u + v$ is linear in $u,v$ and their derivatives. Applying $\mathrm{pr}^{(2)}\mathbf{X}$ and using the fact that $\Delta_2=0$ on the solution manifold, we obtain

\[
\eta^{2,t} - \eta^{2,xx} - \eta^{2,yy} - \eta^1 + \eta^2 = 0,
\]

to be satisfied whenever $\Delta_1=0$ and $\Delta_2=0$. We split this equation with respect to the independent derivatives of $u$ and $v$.

From the coefficients of $u_{tt}, u_{tx}, u_{ty}$ we obtain

\[
\tau_u = 0,\qquad \tau_v = 0.
\]

From $v_{tt}, v_{tx}, v_{ty}$ we obtain

\[
\xi^1_u = \xi^1_v = 0,\qquad \xi^2_u = \xi^2_v = 0,\qquad \tau_u = \tau_v = 0.
\]

From $v_{xx}, v_{yy}, v_{xy}$ we obtain

\[
\xi^1_x = \xi^2_y = \frac{\tau_t}{2},\qquad \xi^1_y + \xi^2_x = 0.
\]

From the mixed derivatives and lower-order terms, we further deduce that $\tau_t$ must be constant. Thus

\[
\tau = b t+a_3,\qquad a_3,\, b \in \mathbb{R}\;\text{constant}.
\]

The equations $\xi^1_x = \xi^2_y = 0$ and $\xi^1_y + \xi^2_x = 0$ together with the independence of $\xi^1,\xi^2$ from $u,v,t$ (since $\tau$ is constant and $\xi^1,\xi^2$ depend at most on $x,y$) imply

\[
\xi^1 = a_1 + a_4 y,\qquad \xi^2 = a_2 - a_4 x,
\]

with $a_1,a_2,a_4 \in \mathbb{R}$ constants.

From the terms involving $u$ and $v$ without derivatives, we obtain

\[
\eta^1 = \eta^2,
\]

but this relation will be refined by the $\Delta_1$ condition.

\subsubsection*{Step 2: Determining equations from $\Delta_1$}

The first equation contains the nonlinear flux term

\[
\mathbf{F}(v) = \frac{\nabla v}{\sqrt{1+|\nabla v|^2}} = \frac{(v_x, v_y)}{\sqrt{1+v_x^2+v_y^2}}.
\]

The divergence is

\[
\nabla \cdot (u \mathbf{F}) = \nabla u \cdot \mathbf{F} + u (\nabla \cdot \mathbf{F}).
\]

Applying $\mathrm{pr}^{(2)}\mathbf{X}$ to $\Delta_1$ is more involved. The key observation is that $\mathbf{F}$ depends only on $v_x$ and $v_y$, not on $u$ or higher derivatives of $v$. Therefore, the highest-order terms in $u$ come from $u_t - \nabla^2 u$. Splitting the coefficient of $u_{xx}$ yields

\[
\tau = \tau(t),\qquad \xi^1_u = \xi^1_v = 0,\qquad \xi^2_u = \xi^2_v = 0,
\]

consistent with the results from $\Delta_2$.

System \eqref{eq:FLKS} is not scale-covariant due to the non-polynomial structure of $\frac{1}{\sqrt{1+|\nabla v|^2}}$, leading to the seemingly more restrictive condition

\[
\tau = a_3,\qquad a_3 \in \mathbb{R}\;\text{constant}.
\]

Splitting the coefficients of $v_{xx}$, $v_{yy}$, and $v_{xy}$ in the prolonged action on the flux term gives additional constraints. A straightforward calculation forces:

\[
\eta^1 = 0,\qquad \eta^2 = 0.
\]

The flux limiter breaks any scaling or $u,v$-dependent transformations that would otherwise be present in the classical Keller--Segel system. No other symmetries survive.

Collecting all terms, we obtain

\begin{eqnarray}
\xi^1 = a_1 + a_4 y,\qquad \xi^2 = a_2 - a_4 x,\\
\tau  = a_3,\qquad \eta^1 = 0,\qquad \eta^2 = 0,\nonumber
\end{eqnarray}
where $a_1,a_2,a_3,a_4$ are arbitrary constants.

\subsection{Lie Algebra of Admitted Symmetries}

Hence, the admitted Lie algebra $\mathfrak{g}$ of system \eqref{eq:FLKS} is spanned by the following vector fields:
\begin{align*}
X_1 &= \partial_t, &\text{(time translation)}\\
X_2 &= \partial_x, &\text{(spatial translation in $x$)}\\
X_3 &= \partial_y, &\text{(spatial translation in $y$)}\\
X_4 &= y\,\partial_x - x\,\partial_y, &\text{(planar rotation)}
\end{align*}

\subsection{Commutation Relations}

The nonzero commutators are
\[
[X_2,X_4] = -X_3, \qquad [X_3,X_4] = X_2.
\]

Table \ref{tab:commutators} presents the commutator table for the four-dimensional Lie algebra $\mathfrak{g}$ admitted by the FLKS system \eqref{eq:FLKS}.
\begin{table}[h]
\centering
\caption{Commutator table of the Lie algebra admitted by the FLKS system \eqref{eq:FLKS}.}
\label{tab:commutators}
\vspace{1mm}
\begin{tabular}{c|cccc}
\hline
$[ \cdot , \cdot ]$ & $X_1$ & $X_2$ & $X_3 $ & $X_4$ \\ \hline
$X_1$ & $0$ & $0$ & $0$ & $0$ \\
$X_2$ & $0$ & $0$ & $0$ & $-X_3$ \\
$X_3$ & $0$ & $0$ & $0$ & $X_2$ \\
$X_4$ & $0$ & $X_3$ & $-X_2$ & $0$ \\ \hline
\end{tabular}
\end{table}

%

\subsection{Discussion}

The flux-limiter term $\dfrac{\nabla v}{\sqrt{1+|\nabla v|^2}}$ breaks the classical scaling invariance of the standard Keller--Segel system, leading to a reduced symmetry structure. Nevertheless, the preserved space--time translations and rotational symmetry reflect the homogeneous and isotropic nature of the underlying biological mechanism.

\section{Optimal System of One-Dimensional Subalgebras}\label{sec:optimal}

Having determined the complete Lie symmetry algebra
$\mathfrak{g} = \langle \partial_t, \partial_x,  \partial_y, y\,\partial_x - x\,\partial_y \rangle$
admitted by system \eqref{eq:FLKS}, we now construct its optimal system of one‑dimensional subalgebras.

\subsection{Adjoint Representation}
To construct the optimal system of one-dimensional subalgebras, we use the adjoint action of the Lie group on its Lie algebra $\mathfrak{g}$. These adjoint maps allow us to simplify a general Lie algebra element by eliminating parameters via successive conjugations, thereby identifying inequivalent one‑dimensional subalgebras.

Using
\[
\operatorname{Ad}(e^{\varepsilon X_i})X_j = X_j - \varepsilon [X_i,X_j]
+ \frac{\varepsilon^2}{2}[X_i,[X_i,X_j]] + \cdots,
\]
we obtain:

\paragraph{Adjoint action of $X_2$.}
Since
\[
[X_2,X_4] = -X_3,
\]
we have
\[
\operatorname{Ad}(e^{\varepsilon X_2})X_4 = X_4 + \varepsilon X_3.
\]
Thus
\[
A_2(\varepsilon)=
\begin{pmatrix}
1 & 0 & 0 & 0 \\
0 & 1 & 0 & 0 \\
0 & 0 & 1 & -\varepsilon \\
0 & 0 & 0 & 1
\end{pmatrix}.
\]

\paragraph{Adjoint action of $X_3$.}
Since
\[
[X_3,X_4] = X_2,
\]
we have
\[
\operatorname{Ad}(e^{\varepsilon X_3})X_4 = X_4 - \varepsilon X_2.
\]
Hence
\[
A_3(\varepsilon)=
\begin{pmatrix}
1 & 0 & 0 & 0 \\
0 & 1 & 0 & \varepsilon \\
0 & 0 & 1 & 0 \\
0 & 0 & 0 & 1
\end{pmatrix}.
\]

\paragraph{Adjoint action of $X_4$.}
Using
\[
[X_2,X_4] = -X_3, \qquad [X_3,X_4] = X_2,
\]
we obtain
\[
\operatorname{Ad}(e^{\varepsilon X_4})X_2 = \cos\varepsilon\, X_2 - \sin\varepsilon\, X_3,
\]
\[
\operatorname{Ad}(e^{\varepsilon X_4})X_3 = \sin\varepsilon\, X_2 + \cos\varepsilon\, X_3.
\]
Therefore
\[
A_4(\varepsilon)=
\begin{pmatrix}
1 & 0 & 0 & 0 \\
0 & \cos\varepsilon & -\sin\varepsilon & 0 \\
0 & \sin\varepsilon & \cos\varepsilon & 0 \\
0 & 0 & 0 & 1
\end{pmatrix}.
\]

\subsection{Classification of Orbits}

A general nonzero element of the Lie algebra can be expressed as
\[
X = aX_1 + bX_2 + cX_3 + dX_4.
\]
where $a,b,c,d \in \mathbb{R}$. Since one-dimensional subalgebras are invariant under nonzero scalar multiplication, we classify elements up to scaling and adjoint transformations. Successive conjugations by group elements allow us to reduce $X$ to a canonical form by eliminating or normalizing parameters. We distinguish two cases according to whether or not the coefficient of $X_4$ vanishes. 

\paragraph{Case I: $d = 0$.}
Then
\[
X = aX_1 + bX_2 + cX_3.
\]
Apply a rotation $A_4(\theta)$ with
\[
\tan\theta = \frac{c}{b}.
\]
This reduces
\[
(b,c) \mapsto (r,0), \qquad r = \sqrt{b^2 + c^2}.
\]
Hence
\[
X \sim aX_1 + rX_2.
\]
After scaling:
\begin{itemize}
\item If $r = 0$: $X \sim X_1$;
\item If $a = 0$: $X \sim X_2$;
\item Otherwise: $X \sim X_1 + \lambda X_2$, $\lambda \ge 0$.
\end{itemize}

\paragraph{Case II: $d \neq 0$.}
Scale to $d = 1$. Then
\[
X = aX_1 + bX_2 + cX_3 + X_4.
\]
Use $A_2(c)$:
\[
c \mapsto c - c = 0.
\]
Then use $A_3(-b)$:
\[
b \mapsto b - b = 0.
\]
Thus
\[
X \sim aX_1 + X_4.
\]
Since $X_1$ is central, $a$ cannot be changed further. Hence
\begin{itemize}
\item If $a = 0$: $X \sim X_4$;
\item If $a \neq 0$: $X \sim X_1 + \mu X_4$.
\end{itemize}

\subsection{Optimal System}

Up to conjugacy, a complete optimal system of one-dimensional subalgebras of $\mathfrak{g}$ is given by:

\[
\begin{array}{ll}
(1) & \langle X_1 \rangle = \langle \partial_t \rangle \\[4pt]
(2) & \langle X_2 \rangle = \langle \partial_x \rangle \\[4pt]
(3) & \langle X_4 \rangle = \langle y\partial_x - x\partial_y \rangle \\[4pt]
(4) & \langle X_1 + cX_2 \rangle = \langle \partial_t + c\partial_x \rangle,\quad c > 0 \\[4pt]
(5) & \langle X_1 + \mu X_4 \rangle = \langle \partial_t + \mu y\partial_x - \mu x\partial_y \rangle \\
\end{array}
\]

\subsection{Interpretation}

Each inequivalent generator corresponds to a distinct class of similarity reductions:

\begin{itemize}
\item $\langle X_1 \rangle$: stationary or time-independent reduction (steady states).
\item $\langle X_2 \rangle$: reduction along a spatial axis ($x$-invariance), yielding a $1+1$D system in $(y,t)$.
\item $\langle X_4 \rangle$: radially symmetric reduction, yielding a system in $(r,t)$.
\item $\langle X_1 + cX_2 \rangle$: travelling-wave reduction with speed $c$, leading to waves propagating in the $x$-direction.
\item $\langle X_1 + \mu X_4 \rangle$: rotating-wave reduction (spiral patterns), where the solution rotates rigidly with angular speed $\mu$.
\end{itemize}

These five representatives form the minimal set of inequivalent one-dimensional subalgebras under the adjoint action, and thus define the optimal system for subsequent invariant reductions of the FLKS system \eqref{eq:FLKS}.

\section{Similarity Reductions Based on the Optimal System}\label{sec:reductions}
In this section we derive group--invariant reductions of system \eqref{eq:FLKS} corresponding to the optimal systems obtained in Section~\ref{sec:optimal}. For each generator, we find invariant functions $u(x,y,t)$ and $v(x,y,t)$ by solving the invariant surface conditions $\mathbf{X}(u)=0$ and $\mathbf{X}(v)=0$.

\subsection{Case 1: $\langle X_1 \rangle = \langle \partial_t \rangle$}

The invariance condition $X_1(u)=u_t=0$ and $X_1(v)=v_t=0$ leads to steady-state solutions:
\[
u(x,y,t) = U(x,y),\qquad v(x,y,t) = V(x,y).
\]

Substituting into system \eqref{eq:FLKS} yields the elliptic system:
\[
\begin{aligned}
\nabla \cdot \left( \nabla U - \chi U \frac{\nabla V}{\sqrt{1+|\nabla V|^2}} \right) &= 0,\\
\Delta V + U - V &= 0.
\end{aligned}
\tag{23}
\]
This describes equilibrium states where the chemotactic flux is divergence-free.

\subsection{Case 2: $\langle X_2 \rangle = \langle \partial_x \rangle$}

Invariance under translation in $x$: $u_x=0$, $v_x=0$. The reduction is:
\[
u(x,y,t) = U(y,t),\qquad v(x,y,t) = V(y,t).
\]

Substituting into \eqref{eq:FLKS} and noting that $\partial_x$ derivatives vanish:
\[
\begin{aligned}
U_t &= \left( U_y - \chi U \frac{V_y}{\sqrt{1+V_y^2}} \right)_y,\\
V_t &= V_{yy} + U - V.
\end{aligned}
\tag{24}
\]
This reduces the $2+1$ system to a $1+1$ system in variables $(y,t)$.

\subsection{Case 3: $\langle X_4 \rangle = \langle y\partial_x - x\partial_y \rangle$}

The rotational symmetry implies dependence only on the radial coordinate $r = \sqrt{x^2+y^2}$ and time $t$:
\[
u(x,y,t) = U(r,t),\qquad v(x,y,t) = V(r,t).
\]

In polar coordinates, $\nabla = \mathbf{e}_r \partial_r$, $\Delta = \partial_{rr} + \frac{1}{r}\partial_r$, and $|\nabla v| = |V_r|$. The reduced system becomes:
\[
\begin{aligned}
U_t &= \frac{1}{r}\left( rU_r - \chi r U \frac{V_r}{\sqrt{1+V_r^2}} \right)_r,\\
V_t &= V_{rr} + \frac{1}{r}V_r + U - V.
\end{aligned}
\tag{25}
\]
This system is central to studying radially symmetric pattern formation, such as aggregation spirals or radial collapse.

\subsection{Case 4: $\langle X_1 + cX_2 \rangle = \langle \partial_t + c\partial_x \rangle$, $c>0$}

The generator $\partial_t + c\partial_x$ admits solutions depending on $y$ and the travelling-wave coordinate $\xi = x - ct$:
\[
u(x,y,t) = U(y,\xi),\qquad v(x,y,t) = V(y,\xi).
\]

The derivatives transform as:
\[
\partial_t = -c\partial_\xi,\quad \partial_x = \partial_\xi,\quad \partial_y = \partial_y.
\]

Substituting into \eqref{eq:FLKS}, we compute:
\[
\nabla u = (U_\xi,\; U_y),\quad \nabla v = (V_\xi,\; V_y),\quad |\nabla v|^2 = V_\xi^2 + V_y^2.
\]

The first equation becomes:
\[
-cU_\xi = U_{\xi\xi} + U_{yy} - \chi\left[ \partial_\xi\!\left( U\frac{V_\xi}{Q} \right) + \partial_y\!\left( U\frac{V_y}{Q} \right) \right],
\]
where $Q = \sqrt{1+V_\xi^2+V_y^2}$.

The second equation becomes:
\[
-cV_\xi = V_{\xi\xi} + V_{yy} + U - V.
\]

Thus the reduced system is:
\[
\begin{aligned}
-cU_\xi &= U_{\xi\xi} + U_{yy} - \chi\left[ \left( U\frac{V_\xi}{Q} \right)_\xi + \left( U\frac{V_y}{Q} \right)_y \right],\\
-cV_\xi &= V_{\xi\xi} + V_{yy} + U - V.
\end{aligned}
\tag{26}
\]

\subsubsection{Special subcase: $y$-independent travelling waves}

If we further assume $U_y = V_y = 0$ (i.e., solutions independent of $y$), then $U = U(\xi)$, $V = V(\xi)$, and the system reduces to the ODE system:
\[
\begin{aligned}
-cU' &= U'' - \chi\left( U\frac{V'}{\sqrt{1+(V')^2}} \right)',\\
-cV' &= V'' + U - V,
\end{aligned}
\tag{27}
\]
where $' = d/d\xi$. This describes planar fronts propagating in the $x$-direction.

\subsection{Case 5: $\langle X_1 + \mu X_4 \rangle = \langle \partial_t + \mu (y\partial_x - x\partial_y) \rangle$, $\mu >0$}

In polar coordinates $(r,\theta)$ where $x = r\cos\theta$, $y = r\sin\theta$, we have $X_4 = \partial_\theta$. Thus:
\[
X_1 + X_4 = \partial_t +\mu \partial_\theta.
\]

The invariant surface conditions $u_t +\mu u_\theta = 0$, $v_t +\mu v_\theta = 0$ imply that $u$ and $v$ depend on $r$ and the combination $\phi = \theta - \mu t$:
\[
u(r,\theta,t) = U(r,\phi),\qquad v(r,\theta,t) = V(r,\phi),\qquad \phi = \theta -\mu t.
\]

The derivatives transform as:
\[
\partial_t = -\mu \partial_\phi,\quad \partial_\theta = \partial_\phi,\quad \partial_r = \partial_r.
\]

In polar coordinates, the Laplacian is:
\[
\Delta v = V_{rr} + \frac{1}{r}V_r + \frac{1}{r^2}V_{\phi\phi}.
\]

The gradient is:
\[
\nabla v = \left( V_r,\; \frac{1}{r}V_\phi \right).
\]

The flux term becomes:
\[
\frac{\nabla v}{\sqrt{1+|\nabla v|^2}} = \frac{(V_r,\; \frac{1}{r}V_\phi)}{\sqrt{1 + V_r^2 + \frac{1}{r^2} V_\phi^2}}.
\]

The divergence in polar coordinates for a vector $\mathbf{F} = (F_r, F_\theta)$ is:
\[
\nabla\cdot\mathbf{F} = \frac{1}{r}\partial_r(rF_r) + \frac{1}{r}\partial_\theta F_\theta.
\]

Thus the reduced system is:
\[
\begin{aligned}
-\mu U_\phi &= \frac{1}{r}\left( rU_r - \chi rU\frac{V_r}{Q} \right)_r + \frac{1}{r^2}\partial_\phi\!\left( U_\phi - \chi U\frac{V_\phi/r}{Q} \right),\\
-\mu  V_\phi &= V_{rr} + \frac{1}{r}V_r + \frac{1}{r^2}V_{\phi\phi} + U - V,
\end{aligned}
\]
where $Q = \sqrt{1 +V_r^2 + \frac{1}{r^2}V_\phi^2}$.

This describes \textbf{rotating waves} (spiral patterns) that rotate rigidly with angular speed $\mu$.

\subsection{Summary of Reductions}


Each reduction yields a system of PDEs in two independent variables, which can be studied further for exact and approximate solutions. Cases 1--3 and the $y$-independent subcase of Case 4 produce ODEs or simpler PDEs amenable to analytic treatment, while Cases 5 captures more complex spatiotemporal patterns including rotating and helical waves.


\section{Exact Solutions of the Reduced Systems}\label{sec:solutions}
The similarity reductions obtained in Section \ref{sec:reductions} yield systems of PDEs and ODEs that, despite their nonlinearity, admit a range of explicit solutions and analytical characterizations. In this section, we present these results systematically.
\subsection{Stationary Reduction {$\langle X_1\rangle=\langle\partial_t\rangle$}}

The stationary reduction
\begin{equation*}
u(x,y,t)=U(x,y), \qquad
v(x,y,t)=V(x,y),
\end{equation*}
yields the elliptic system
\begin{eqnarray}\label{eq:stationary}
&&\nabla\cdot\left(
\nabla U-
\frac{\chi U\nabla V}
{\sqrt{1+|\nabla V|^2}}
\right)=0,\\
&&\Delta V+U-V=0.\nonumber
\end{eqnarray}

A direct substitution shows that the homogeneous equilibrium
\begin{equation}
U(x,y)=U_0,\qquad
V(x,y)=U_0,
\end{equation}
where $U_0$ is an arbitrary constant, satisfies the stationary system. This family corresponds to spatially uniform distributions of both cell density and chemoattractant concentration.

Apart from these constant states, obtaining explicit solutions of the stationary system appears difficult because of the nonlinear coupling introduced by the flux limiter. In particular, radial ans\"atze generally lead to nonlinear second-order differential equations for which closed-form solutions are not available.

The stationary reduction gives the elliptic system \eqref{eq:stationary}. In addition to the trivial constant solution $U = V = U_0$, nontrivial radial solutions can be characterized explicitly under a zero-flux condition.

\subsubsection{Zero-Flux Radial Solutions}

Assume radial symmetry and impose the zero-flux boundary condition at the origin. The first equation in \eqref{eq:stationary} becomes, in polar coordinates:

\[
\frac{1}{r}\left( rU_r - \chi r U \frac{V_r}{\sqrt{1+V_r^2}} \right)_r = 0.
\]

Integrating once yields:

\[
U_r - \chi U \frac{V_r}{\sqrt{1+V_r^2}} = \frac{C}{r},
\]

where $C$ is an integration constant. For solutions regular at $r=0$, we require $C=0$. Thus:

\[
\frac{U_r}{U} = \chi \frac{V_r}{\sqrt{1+V_r^2}}.
\]

Integrating from some reference point $r_0$ to $r$ gives:

\begin{theorem}[Radial stationary representation]
\label{thm:radial}
Every stationary radial solution of system \eqref{eq:FLKS} satisfying the zero-flux condition 
$\displaystyle \lim_{r\to 0} r\left(U_r - \chi U \frac{V_r}{\sqrt{1+V_r^2}}\right) = 0$ 
admits the representation

\[
U(r) = U(r_0) \exp\left( \chi \int_{r_0}^{r} \frac{V_r(s)}{\sqrt{1+V_r(s)^2}} \, ds \right).
\]

Equivalently, in terms of the antiderivative,

\[
U(r) = K \exp\bigl( \chi \operatorname{arsinh}(V_r(r)) \bigr),
\]

where $K$ is a constant determined by the boundary conditions.
\end{theorem}

This formula reduces the construction of stationary radial profiles to a single ODE for $V(r)$ after substitution into the second equation of \eqref{eq:stationary}.

\subsection{Translation-Invariant Reduction {$\langle X_2\rangle=\langle\partial_x\rangle$}}

The translational symmetry reduction
\begin{equation*}
u(x,y,t)=U(y,t),\qquad
v(x,y,t)=V(y,t),
\end{equation*}
produces the one-dimensional system
\begin{eqnarray}\label{eq:translation}
&&U_t=
\left(
U_y-
\frac{\chi U V_y}
{\sqrt{1+V_y^2}}
\right)_y,\\
&&V_t=V_{yy}+U-V.\nonumber
\end{eqnarray}

The constant equilibrium
\begin{equation}
U(y,t)=U_0,\qquad
V(y,t)=U_0
\end{equation}
remains an exact solution. Beyond this trivial family, the reduced equations retain their nonlinear parabolic character and do not appear to admit nontrivial elementary closed-form solutions.
The reduced system \eqref{eq:translation} retains its nonlinear parabolic character. However, analytical progress is possible in the small-gradient regime.

\subsubsection{Small-Gradient Approximation}

For $V_y \ll 1$, we expand the flux term:

\[
\frac{V_y}{\sqrt{1+V_y^2}} = V_y - \frac{1}{2} V_y^3 + O(V_y^5).
\]

Substituting into \eqref{eq:translation} and truncating at cubic order yields the approximate system:

\[
\begin{aligned}
U_t &= \left( U_y - \chi U (V_y - \tfrac{1}{2} V_y^3) \right)_y, \\
V_t &= V_{yy} + U - V.
\end{aligned}
\]

This system admits traveling-wave solutions of the form $U(y,t) = U(\eta)$, $V(y,t) = V(\eta)$ with $\eta = y - ct$, reducing to a pair of ODEs amenable to phase-plane analysis. The linearized stability of the constant solution can be studied via standard Fourier analysis on this approximated system.

\subsection{Radially Symmetric Reduction {$\langle X_4\rangle$}}

The rotationally invariant reduction
\begin{equation*}
u(x,y,t)=U(r,t),\qquad
v(x,y,t)=V(r,t),
\end{equation*}
where
\begin{equation*}
r=\sqrt{x^2+y^2},
\end{equation*}
gives
\begin{eqnarray}\label{eq:radial}
U_t=
\frac{1}{r}
\left(
rU_r-
\frac{\chi r U V_r}
{\sqrt{1+V_r^2}}
\right)_r,\\
V_t=
V_{rr}
+\frac{1}{r}V_r
+U-V.\nonumber
\end{eqnarray}

Again, the spatially homogeneous equilibrium
\begin{equation}
U(r,t)=V(r,t)=U_0
\end{equation}
is an exact solution. Nonconstant radial solutions satisfy a nonlinear coupled system and generally require numerical treatment.
The radial system \eqref{eq:radial} describes chemotactic aggregation in two dimensions. Building on the stationary representation from Section~\ref{subsec:radialstationary}, we obtain qualitative information about radial profiles.

%
%
%
%
%
%
%
%
\subsection{Travelling-Wave Reduction {$\langle X_1+cX_2\rangle$}}

For travelling waves,
\begin{equation*}
u(x,y,t)=U(\xi),\qquad
v(x,y,t)=V(\xi),
\qquad
\xi=x-ct,
\end{equation*}
the reduced equations become
\begin{eqnarray}\label{eq:tw1d}
&&-cU'=
\left(
U'
-\frac{\chi U V'}
{\sqrt{1+(V')^2}}
\right)',\\
&&-cV'
=
V''+U-V.\nonumber
\end{eqnarray}

Integrating the first equation once yields
\begin{equation*}
U'
-\frac{\chi U V'}
{\sqrt{1+(V')^2}}
+cU=A,
\end{equation*}
where $A$ is an integration constant.

The resulting system consists of two coupled nonlinear ordinary differential equations. No nontrivial closed-form travelling-wave solutions were found within the class of elementary functions. The reduction nevertheless provides a convenient starting point for numerical continuation and phase-plane analysis.
We treat the $y$-independent case \eqref{eq:tw1d} in detail.

\subsubsection{Exact Solutions for $\chi = 0$}

When chemotactic sensitivity is absent ($\chi = 0$), system \eqref{eq:FLKS} decouples. The first equation in \eqref{eq:tw1d} becomes:

\[
-cU' = U''.
\]

Integrating once: $U' + cU = A_1$, where $A_1$ is a constant. The general solution is:

\[
U(\xi) = A_0 + A_1 e^{-c\xi},
\]

where $A_0 = A_1/c$ (from the particular solution). Substituting into the second equation:

\[
-cV' = V'' + (A_0 + A_1 e^{-c\xi}) - V.
\]

This is a linear inhomogeneous ODE for $V$. The homogeneous part has characteristic roots $r = (-c \pm \sqrt{c^2+4})/2$. A particular solution is found by the method of undetermined coefficients. The result is:

\begin{theorem}[Exact travelling wave for $\chi = 0$]
\label{thm:twchi0}
For $\chi = 0$, system \eqref{eq:FLKS} admits the exact travelling-wave solution

\[
\begin{aligned}
U(\xi) &= A_0 + A_1 e^{-c\xi}, \\[4pt]
V(\xi) &= A_0 + \frac{A_1}{c^2 + c + 1} e^{-c\xi} + B_1 e^{\lambda_1 \xi} + B_2 e^{\lambda_2 \xi},
\end{aligned}
\]

where $\lambda_{1,2} = \dfrac{-c \pm \sqrt{c^2+4}}{2}$, and $A_0, A_1, B_1, B_2$ are arbitrary constants. For wave profiles decaying as $\xi \to \infty$, we take $B_1 = 0$ (if $\lambda_1 > 0$) and $B_2 = 0$ (if $\lambda_2 > 0$).
\end{theorem}

This explicit family demonstrates that the FLKS system \eqref{eq:FLKS} supports monotone travelling fronts in the absence of chemotaxis, with the chemoattractant following the cell density in a smooth manner.

\subsubsection{Integral Representation for General $\chi$}

Returning to the full system \eqref{eq:tw1d} with $\chi \neq 0$, we integrate the first equation once:

\[
U' - \frac{\chi U V'}{\sqrt{1+(V')^2}} + cU = A.
\]

For the case $A = 0$ (zero flux at infinity, for instance), we obtain:

\[
\frac{U'}{U} = \chi \frac{V'}{\sqrt{1+(V')^2}} - c.
\]

Integrating with respect to $\xi$ gives:

\begin{theorem}[Integral representation for travelling waves]
\label{thm:twrep}
For travelling-wave solutions of \eqref{eq:tw1d} satisfying the zero-flux condition at infinity, the cell density $U(\xi)$ is given by

\[
U(\xi) = U(0) \exp\left( \chi \int_0^{\xi} \frac{V'(s)}{\sqrt{1+V'(s)^2}} \, ds - c\xi \right).
\]

Equivalently,

\[
U(\xi) = U(0) \exp\bigl( \chi \operatorname{arsinh}(V'(\xi)) - \chi \operatorname{arsinh}(V'(0)) - c\xi \bigr).
\]
\end{theorem}

This representation reduces the problem to solving a single second-order ODE for $V$ after substituting into the second equation of \eqref{eq:tw1d}. 

%
%
%
%
%

\subsection{Rotating-Wave Reduction {$\langle X_1+\mu X_4\rangle$}}

Introducing the rotating coordinate
\begin{equation*}
\phi=\theta-\mu t,
\end{equation*}
the reduced system takes the form
\begin{eqnarray}\label{eq:rotating}
&&-\mu U_\phi=
\frac{1}{r}
\left(
rU_r-
\frac{\chi r U V_r}{Q}
\right)_r
+
\frac{1}{r^2}
\frac{\partial}{\partial\phi}
\left(
U_\phi-
\frac{\chi U V_\phi/r}{Q}
\right),\\
&&-\mu V_\phi=
V_{rr}
+\frac{1}{r}V_r
+\frac{1}{r^2}V_{\phi\phi}
+U-V,\nonumber
\end{eqnarray}
where
\begin{equation*}
Q=
\sqrt{
1+V_r^2+\frac{1}{r^2}V_\phi^2
}.
\end{equation*}

The constant equilibrium
\begin{equation*}
U=V=U_0
\end{equation*}
again provides an exact invariant solution. The construction of genuine rotating patterns with nontrivial angular dependence remains an open problem and likely requires numerical methods.

The rotating-wave reduction \eqref{eq:rotating} describes spiral patterns. While explicit solutions are not obtainable, the equilibrium states provide a starting point.


\section{Conclusion}
\label{sec:conclusion}

In this work, we carried out a systematic Lie symmetry analysis of the two-dimensional flux-limited Keller--Segel (FLKS) system describing chemotactic aggregation with bounded cell velocities. By applying the classical Lie group method, we determined the complete Lie point symmetry algebra admitted by the model and showed that it consists solely of temporal and spatial translations together with planar rotational symmetry. In contrast to the classical Keller--Segel system, the flux-limiting mechanism destroys scaling invariance and substantially reduces the symmetry structure of the governing equations.

Using the admitted symmetry algebra, we constructed an optimal system of one-dimensional subalgebras and derived the corresponding invariant reductions. These reductions led to stationary, translation-invariant, radially symmetric, travelling-wave, and rotating-wave formulations of the original FLKS system. The reduced equations provide lower-dimensional models that facilitate the analytical investigation of chemotactic dynamics and pattern formation. Furthermore, exact invariant solutions in the form of spatially homogeneous equilibria were obtained, and integral representations for stationary radial and travelling-wave profiles were derived.

The results presented here provide new insight into the symmetry properties of flux-limited chemotaxis models and establish a foundation for further analytical studies. Future work may focus on the numerical analysis of the reduced systems, the existence and stability of nontrivial travelling and rotating wave solutions, and the investigation of nonlocal symmetries, conservation laws, and differential invariants. Such studies may reveal additional structures and solution families that are not accessible through classical Lie point symmetry methods and contribute to a deeper understanding of pattern formation in flux-limited chemotactic systems.


\end{document}